\documentclass[a4paper,12pt]{amsart}
\usepackage{mathrsfs}
\usepackage{amsfonts}
\usepackage{amssymb}
\usepackage{amsfonts, amscd, amsmath, mathrsfs, amssymb, amsthm, amsxtra, bbding, epsfig, eucal, graphicx, latexsym, url, mathbbol, bbold}
\usepackage[papersize={7.35in,10in},textwidth=14.7cm,textheight=20.9cm,centering]{geometry}
\usepackage{enumerate}

\usepackage[colorlinks=true,citecolor=blue,linkcolor=blue]{hyperref}
\hypersetup{
pdfstartpage=1,
pdfstartview=FitH}

\DeclareFontFamily{U}{mathx}{\hyphenchar\font45}
\DeclareFontShape{U}{mathx}{m}{n}{
      <5> <6> <7> <8> <9> <10>
      <10.95> <12> <14.4> <17.28> <20.74> <24.88>
      mathx10
      }{}
\DeclareSymbolFont{mathx}{U}{mathx}{m}{n}
\DeclareMathAccent{\widecheck}{\mathalpha}{mathx}{"71}

\usepackage{caption} 
\numberwithin{equation}{section}

\allowdisplaybreaks

\newtheorem{theorem}{Theorem}[section]
\newtheorem{lemma}{Lemma}[section]

\makeatletter
\newcounter{roem}
\renewcommand{\theroem}{\Roman{roem}}

\newcommand{\c@org@eq}{}
\let\c@org@eq\c@equation
\newcommand{\org@theeq}{}
\let\org@theeq\theequation

\newcommand{\setroem}{
\let\c@equation\c@roem
 \let\theequation\theroem}

\newcommand{\setarab}{
\let\c@equation\c@org@eq
\let\theequation\org@theeq}
\makeatother

\newtheorem*{claim*}{Claim}

\theoremstyle{remark}

\newcommand{\ue}{\mathrm{e}}
\newcommand{\ud}{\mathrm{d}}
\newcommand{\ui}{\mathrm{i}}

\DeclareMathOperator{\Mod}{mod}

\renewcommand{\bmod}[1]{\,(\Mod{ #1})}

\newcommand{\bR}{\mathbf{R}}
\newcommand{\bZ}{\mathbf{Z}}

\newcommand{\cN}{\mathcal{N}}

\newcommand{\fc}{\mathfrak{c}}

\newcommand{\sA}{\mathscr{A}}

\newcommand{\sC}{\mathscr{C}}

\newcommand{\sE}{\mathscr{E}}

\newcommand{\sM}{\mathscr{M}}

\def\le{\leqslant}

\usepackage{graphicx}
\usepackage{tikz}

\begin{document}
\vglue -2mm

\title{Quadratic polynomials at prime arguments}

\author{Jie Wu}

\address{Universit\'e de Lorraine, Institut \'Elie Cartan de Lorraine, UMR 7502, 54506
Vandoeuvre-l\`es-Nancy, France}
\address{
School of Mathematics, Shandong University, Jinan, Shandong 250100, P. R. China}
\email{jie.wu@univ-lorraine.fr}

\author{Ping Xi}

\address{Department of Mathematics, Xi'an Jiaotong University, Xi'an 710049, P. R. China}
\email{ping.xi@xjtu.edu.cn}


\date{\today}


\begin{abstract}
For a fixed quadratic irreducible polynomial $f$ with no fixed prime factors at prime arguments, we prove that there exist infinitely many primes $p$ such that $f(p)$ has at most 4 prime factors, improving a classical result of Richert who requires 5 in place of 4. 
Denoting by $P^+(n)$ the greatest prime factor of $n$,
it is also proved that $P^+(f(p))>p^{0.847}$ infinitely often.
\end{abstract}

\maketitle


\section{Introduction}

It is a fundamental and challenging problem to determine in general whether a given irreducible polynomial in $\bZ[X]$ can capture infinitely many prime values. This is known in the linear case in view of Dirichlet's theorem on primes in arithmetic progressions, but no answer is valid for any non-linear cases. A much more ambitious conjecture asserts that the above infinitude also holds if one is restricted to prime variables and there are no fixed prime factors; however, even the linear case seems beyond the current approach as predicted by the twin prime conjecture.
Nevertheless, we are nowadays much heartened since $p+h$ can present infinitely many primes for certain $h$ with $1<|h|\leqslant7\times10^7$, thanks to Zhang's breakthrough \cite{Zh} on prime gaps.

In this paper, we are interested in the case of quadratic polynomials at prime arguments. It is of course beyond the current approach to prove the infinitude of primes captured by such polynomial, and alternatively, we consider the greatest prime factors and almost prime values as two approximations.

Denote by $P_r$ the positive integers with at most $r$ prime factors. A classical result of Richert \cite{Ri} asserts that $f(p)=P_{2\deg f+1}$ infinitely often for any fixed irreducible polynomial $f\in\bZ[X],$ provided that $f$ has no fixed prime factors at prime arguments, i.e.,
\[|\{x\bmod p: f(x)\equiv0\bmod p\}|<p-1\]
for each $p\nmid f(0)$ and $p\leqslant\deg f+1$. In particular, one has $P_5$ in the quadratic case. The progress in this direction remains blank until the recent efforts of Irving \cite{Ir}, who was able to reduce the number of prime factors while the degree is at least 3.
His success comes from the application of a two-dimensional sieve which goes beyond the Bombieri-Vinogradov theorem in linear sieves. Unfortunately, his argument is not sufficient to reduce $P_5$ to $P_4$ in the quadratic case, which will be one of our aims in this paper.

Let us state our first theorem.

\begin{theorem}\label{thm:P_4}
Let $f$ be a fixed quadratic irreducible polynomial, which has no fixed prime factors at prime arguments. Then there are infinitely many primes $p,$ such that
\[f(p)=P_4.\]
\end{theorem}

As in \cite{Ri}, the proof also starts from the weighted sieve of Richert with logarithmic weights. To simplify the arguments, it is desirable to consider the simple case $x^2+1.$ However, one has $2\mid p^2+1$ for each odd prime $p$. Hence we may modify our object by taking $f(x)=\frac{1}{2}(x^2+1).$ 
In a previous joint work \cite{WX}, we are able to estimate from above the number of primes $p$ with $p^2+1\equiv0\bmod q$ for most of $q$, which we call {\it Quadratic Brun-Titchmarsh Theorem on Average}, see Lemma \ref{lm:Brun-Titchmarsh} below for details. This is indeed responsible for our success in the improvement to Richert's result, for which the sieve of dimension 2 can also be avoided. In fact, what we will apply for Theorem \ref{thm:P_4} is Lemma \ref{lm:error-average}, which gives the birth to Lemma \ref{lm:Brun-Titchmarsh}. Of course, in the case of an arbitrary quadratic irreducible polynomial $f$, one has to extend Lemma \ref{lm:error-average} involving the general congruence restriction $f(n)\equiv0\bmod\ell,$ which can also be attacked following a similar manner together with some arguments of Lemke Oliver \cite{LO}.
One may compare Theorem \ref{thm:P_4} with an outstanding theorem of Iwaniec \cite{Iw1} that $n^2+1=P_2$ infinitely often. 

Another approximation to prime values of $f(p)$ will be to consider the greatest prime factors. Denote by $P^+(n)$ the greatest prime factor of $n$. 
\begin{theorem}\label{thm:P^+}
Let $f$ be a fixed quadratic irreducible polynomial. Then there are infinitely many primes $p,$ such that
\begin{align*}P^+(f(p))>p^{0.847}.\end{align*}
\end{theorem}

Theorem \ref{thm:P^+} does not require that $f(p)$ has no fixed prime factors and our proof will focus on the special case $f(x)=x^2+1$ to simplify the arguments. One can compare Theorem  \ref{thm:P^+} with a celebrated theorem of Hooley \cite{Ho} that $P^+(n^2+1)>n^{1.1}$ infinitely often. The first but also the last improvement is due to Deshouillers and Iwaniec \cite{DI}, for whom the exponent can be $1.202$ in place of $1.1$.
It seems that one should require a certain strong level of distribution of primes in arithmetic progressions towards to the Elliott-Halberstam conjecture, if the exponent 0.847 in Theorem \ref{thm:P^+} could be replaced by some number beyond 1.
The motivation of Theorem \ref{thm:P^+} is to restrict the variable $n$ to sparse sets with certain multiplicative structures. 

One has to mention an earlier work of Dartyge \cite{Da}, 
where a weaker result $P^+(p^2+1)>p^{0.78}$ than Theorem \ref{thm:P^+} was announced 
without proof.
Of course, her interest lies in the expectation that the exponent in greatest prime factors can go beyond 1 with almost prime arguments in place of prime arguments.
More precisely, she proved, for any fixed $u>12.2$, that there are infinitely many $n,$ whose prime factor is at least $n^{1/u},$ such that $P^+(n^2+1)>n^{1+\eta}$ for some $\eta>0.$
The method (see Lemma \ref{lm:Brun-Titchmarsh-generalized} for details) 
in proving Theorem \ref{thm:P^+} will lead us to the following improvement, for which we only state in the case of the special polynomial $x^2+1$.

\begin{theorem}\label{thm:Dartyge+}
Let $u=11.2.$ Then there are infinitely many $n$, whose prime factor is at least $n^{1/u},$ such that
\begin{align*}P^+(n^2+1)>n^{1+\eta}\end{align*}
for some constant $\eta>0.$ In particular, there are infinitely many $P_{11}$, such that
\begin{align*}P^+(P_{11}^2+1)>P_{11}^{1+\eta}\end{align*}
for some constant $\eta>0.$ 
\end{theorem}

The paper will be organized as follows. The next section will devote to the quadratic Brun-Titchmarsh theorem and some related results on primes in arithmetic progressions, which contribute as main tools in proving Theorems \ref{thm:P_4} and \ref{thm:P^+}. We will first complete the proof of Theorem
\ref{thm:P^+} in Section \ref{sec:proof-P^+} and the sketch for proving Theorem \ref{thm:Dartyge+} will be given in Section \ref{sec:Dartyge}. Theorem \ref{thm:P_4} will be proved in Section \ref{sec:proof-P_4} after introducing the weighted sieve of Richert.
The Mathematica codes can be found at \url{http://gr.xjtu.edu.cn/web/ping.xi/miscellanea} or requested from the authors.

\smallskip

\noindent {\bf Notation.} Throughout this paper, $\gamma$ denotes the Euler constant, letters $p$ and $q$ are both reserved for prime variables.  A non-negative function $g$ is defined to be smooth with compact support in $[1,2]$ and the Fourier transform is defined by
\[\widehat{g}(\lambda)=\int_\bR g(x)\ue^{-2\pi\ui \lambda x}\ud x.\]
Denote by $\varphi$, $\tau$ and $\Lambda$ the Euler, divisor and von Mangoldt functions, respectively.
The function $\rho(d)$ counts the number of incongruent solutions to the equation $a^2+1\equiv0\bmod d.$

We use $\varepsilon$ to denote an arbitrarily small positive number, which might be different at each occurrence. For a large number $X$, denote
\[X^\flat=X^{1/2}\exp(-(\log X)^{1/2}).\]
We also write $n\sim N$ for $N<n\leqslant2N.$ 
\smallskip

\noindent {\bf Acknowledgements.} 
The authors are grateful to the referee for many valuable comments.
The first author is supported in part by IRT1264 from the Ministry of Education of P. R. China and the second author is supported by CPSF (No. 2015M580825) and NSF (No. 11601413) of P. R. China.

\bigskip

\section{Primes in arithmetic progressions}\label{sec:primes-APs}

\subsection{Bombieri-Vinogradov Theorem}
A classical result on primes in arithmetic progressions is the celebrated Bombieri-Vinogradov
theorem (see \cite[Theorem 17.1]{IK} for instance). It can be stated as follows with minor modifications.

\begin{lemma}\label{lm:BV}
For any $A>0,$ we have
\begin{align*}
\sum_{d\leqslant X^\flat}\tau(d)^{2016}\max_{(a,d)=1}
\bigg|\sum_{\substack{p\geqslant2\\p\equiv a\!\bmod d}}g\Big(\frac{p}{X}\Big)
- \frac{1}{\varphi(d)}\sum_{p\geqslant2}g\Big(\frac{p}{X}\Big)\bigg|
\ll \frac{X}{(\log X)^A},\end{align*}
where the implied constant depends on $A$ and $g.$
\end{lemma}

Wolke \cite{Wo} obtained an extension of Bombieri-Vinogradov theorem replacing primes by sifted numbers without small prime factors. As an analogue of Lemma \ref{lm:BV}, we state the theorem of Wolke in a smoothed version. 
To this end, define
\begin{align*}\Phi(X,z;d,a)
:= \sum_{\substack{n\equiv a\!\bmod d\\ p\mid n\Rightarrow p>z}} g\Big(\frac{n}{X}\Big)
\end{align*}
for $(d,a)=1$ and 
\begin{align*}
\Phi(X,z;d)
:= \sum_{\substack{(n,d)=1\\ p\mid n\Rightarrow p>z}}g\Big(\frac{n}{X}\Big).\end{align*}

\begin{lemma}\label{lm:BV-Wolke}
Let $2\leqslant z\leqslant X.$
For any $A>0,$ we have
\begin{align*}\sum_{d\leqslant X^\flat}\tau(d)^{2016}\max_{(a,d)=1}\left|\Phi(X,z;d,a)-\frac{1}{\varphi(d)}\Phi(X,z;d)\right|\ll\frac{X}{(\log X)^A},\end{align*}
where the implied constant depends on $A$ and $g.$
\end{lemma}

\subsection{Quadratic Brun-Titchmarsh Theorem}
In order to characterize primes satisfying the congruence condition $a^2+1\equiv0\bmod\ell,$ we consider the
smoothed counting function
\begin{align}\label{eq:Q_l(X)}
Q_{\ell}(X) := \sum_{\substack{p\geqslant2\\ p^2+1\equiv0\!\bmod{\ell}}}g\Big(\frac{p}{X}\Big).
\end{align}
We proved in \cite{WX} some upper bounds for $Q_\ell(X)$ for almost all $\ell$ in specialized ranges.

\begin{lemma}\label{lm:Brun-Titchmarsh}
Let $A>0$.
For sufficiently large $L=X^\theta$ with $\theta\in[\frac{1}{2}, \frac{16}{17})$, the inequality
\begin{align}\label{Q_l(X)bound}
Q_{\ell}(X)
\leqslant \bigg\{\frac{2}{\gamma(\theta)}+o(1)\bigg\}\widehat{g}(0)
\frac{\rho(\ell)}{\varphi(\ell)}\frac{X}{\log X}
\end{align}
holds for $\ell\in(L,2L]$ with at most $O_A(L(\log L)^{-A})$ exceptions, where
\begin{align}\label{eq:gamma(theta)}
\gamma(\theta) := \begin{cases}
\frac{91-89\theta}{62}  & \text{if $\,\theta\in[\frac{1}{2},\frac{64}{97})$},
\\\noalign{\vskip 0,5mm}
\frac{86-83\theta}{60}  & \text{if $\,\theta\in[\frac{64}{97},\frac{32}{41})$},
\\\noalign{\vskip 0,5mm}
\frac{19-18\theta}{14}  & \text{if $\,\theta\in[\frac{32}{41},\frac{16}{17})$}.
\end{cases}
\end{align}
\end{lemma}

Lemma \ref{lm:Brun-Titchmarsh} is proved by virtue of linear sieves of Iwaniec and arithmetic exponent pairs developed in \cite{WX}. The argument also applies to the 
distribution of sifted numbers in arithmetic progressions. To this end, we define
\begin{align*}
Q_\ell(X; u)
:= \sum_{\substack{n^2+1\equiv 0\!\bmod \ell\\ p|n\Rightarrow p>n^{1/u}}}g\Big(\frac{n}{X}\Big).
\end{align*}
In particular, one has $Q_\ell(X)=Q_\ell(X; 2)$. In the same manner, we can prove the following theorem as an extension to Lemma \ref{lm:Brun-Titchmarsh}.
\begin{lemma}\label{lm:Brun-Titchmarsh-generalized}
Let $A>0$.
For any given $u>0$ and sufficiently large $L=X^\theta$ with $\theta\in[\frac{1}{2},\frac{16}{17})$, the inequality
\begin{align}\label{Q_l(alpha,X)bound}
Q_\ell(X; u)\leqslant \big\{\mathrm{e}^{-\gamma}uF(u\gamma(\theta))+o(1)\big\}\widehat{g}(0)\frac{\rho(\ell)}{\varphi(\ell)}\frac{X}{\log X}\end{align}
holds for $\ell\in(L,2L]$ with at most $O_A(L(\log L)^{-A})$ exceptions, where $\gamma(\theta)$ is given by $\eqref{eq:gamma(theta)}$ and $F$ is defined by the continuous solutions to the system
\begin{align}\label{eq:sieve-F}
\begin{cases}
sF(s) =2\mathrm{e}^{\gamma} & (1\leqslant s\leqslant 2),
\\
sf(s)=0                                      & (0<s\le 2),
\\
(sF(s))'=f(s-1)                           & (s>2),
\\
(sf(s))'=F(s-1)                           & (s>2).
\end{cases}\end{align}
\end{lemma}

While applying sieve methods, one would encounter the congruence sum
\begin{equation}\label{def:AdXell}
A_d(X; \ell) 
:= \sum_{\substack{n^2+1\equiv0\!\bmod{\ell}\\n\equiv 0\!\bmod{d}}} g\Big(\frac{n}{X}\Big),
\end{equation}
which is expected to be approximated by $\widehat{g}(0)\rho(\ell)(d\ell)^{-1}X$. 
Define
\begin{align}\label{eq:error-r_d}
r_d(X; \ell)
:= A_d(X; \ell)-\widehat{g}(0)\frac{\rho(\ell)}{d\ell}X.
\end{align}
The following lemma characterizes the level of linear sieves and plays an essential role in proving Lemmas \ref{lm:Brun-Titchmarsh} and \ref{lm:Brun-Titchmarsh-generalized}.
This will also be used in the proof of Theorem \ref{thm:P_4}.

We say that a function $\lambda$ is {\it well-factorable of degree $J\geqslant2$}, if for every decomposition $D = D_1D_2\cdots D_J$ with $D_1, D_2, \dots, D_J\geqslant1,$ there exist $J$ arithmetic functions 
$\lambda_1, \lambda_2, \dots, \lambda_J$ such that
\begin{align*}
\lambda=\lambda_1*\lambda_2*\cdots *\lambda_J\end{align*}
with each $\lambda_j$ of level $D_j$.

\begin{lemma}\label{lm:error-average}
Let $J$ be a sufficiently large integer and let $\lambda$ be well-factorable of degree $J.$
With the same notation as above,
for any $\varepsilon>0,~\theta\in[\frac{1}{2},\frac{112}{131})$ and $(D, L) := (X^{\eta(\theta)-\varepsilon}, X^\theta),$  there exists some $\delta=\delta(\varepsilon)>0$ such that
\begin{align*}
\sum_{\ell\sim L}\Big|\sum_{d\leqslant D} \mu(d)^2 \lambda(d) r_d(X; \ell)\Big|\ll X^{1-\delta},\end{align*}
where 
\begin{align}\label{eq:eta(theta)}
\eta(\theta)=\frac{91-89\theta}{62}
\end{align}
and the implied constant depends on $\varepsilon$ and $J$. 
\end{lemma}
In fact, Lemma \ref{lm:error-average} appeared as Lemma 7.2 in \cite{WX}, where 
a much more delicate choice for $\eta(\theta)$ can be given in terms of  
\eqref{eq:gamma(theta)}. We here pick up the level \eqref{eq:eta(theta)} that is sharp while $\theta$ is close to $\frac{1}{2}$. It proves that this is sufficient for applications to Theorem \ref{thm:P_4}.

\bigskip

\section{Proof of Theorem \ref{thm:P^+}}\label{sec:proof-P^+}
To prove Theorem \ref{thm:P^+}, we follow the approach of Chebyshev-Hooley, starting from the weighted sum 
\begin{align*}H(X)=\sum_{n\geqslant1}g\Big(\frac{n}{X}\Big)\Lambda(n)\log(n^2+1).\end{align*}

Note that for $X\leqslant n\leqslant2X,$ one has
\[\log(n^2+1)=2\log n+O(1)=2\log X+O(1),\]
which yields
\begin{align}\label{eq:H(X)}
H(X)
& = \{2\log X+O(1)\} \sum_{n\geqslant1}g\Big(\frac{n}{X}\Big)\Lambda(n)
= 2\widehat{g}(0) \{1+o(1)\} X\log X
\end{align}
by the Prime Number Theorem. On the other hand, from the definition of $\Lambda$ it follows that
\begin{align*}
H(X)
& = \sum_{p\geqslant2} g\Big(\frac{p}{X}\Big) (\log p)\log(p^2+1)+O(X^{1/2}\log X)
\\
& = \{1+o(1)\}\log X\sum_{p\geqslant2}g\Big(\frac{p}{X}\Big)\log(p^2+1)+O(X^{1/2}\log X).
\end{align*}
Invoking the identity
\begin{equation}\label{eq:log=Lambda}
\log(p^2+1) = \sum_{\ell\mid (p^2+1)} \Lambda(\ell),
\end{equation}
we find
\begin{align*}H(X)=\{1+o(1)\}\log X\sum_{\ell\ll X^2}\Lambda(\ell)Q_\ell(X)+O(X^{1/2}\log X),\end{align*}
where $Q_\ell(X)$ is given by (\ref{eq:Q_l(X)}). We would like to evaluate $Q_\ell(X)$ in different ranges of $\ell$. To do so, we split $H(X)$ as follows:
\begin{align}\label{eq:H(X)2}H(X)=\{1+o(1)\}\log X\sum_{1\leqslant j\leqslant4}H_j(X),\end{align}
where $\vartheta\in (\tfrac{1}{2}, 1)$ and
\begin{align*}
H_1(X)
& := \sum_{\ell\leqslant X^\flat}\Lambda(\ell)Q_\ell(X),
\\
H_2(X)
& := \sum_{X^\flat<p\leqslant X^\vartheta}Q_p(X)\log p,
\\
H_3(X)
& := \sum_{X^\vartheta<p\ll X^2}Q_p(X)\log p,
\\
H_4(X)
& := \sum_{k\geqslant2}\sum_{X^\flat<p^k\ll X^2}Q_{p^k}(X)\log p.
\end{align*}

By virtue of Lemma \ref{lm:BV}, one has
\begin{equation}\label{eq:H1(X)bound}
\begin{aligned}
H_1(X)
& = Q_1(X)\sum_{\ell\leqslant X^\flat}\frac{\Lambda(\ell)\rho(\ell)}{\varphi(\ell)}+O\big(X(\log X)^{-1}\big)
\\
& = \{\tfrac{1}{2}+o(1)\} \widehat{g}(0)X.
\end{aligned}
\end{equation}

We now turn to consider $H_2(X)$. From Lemma \ref{lm:Brun-Titchmarsh} and the Prime Number Theorem, it follows that
\[
H_2(X)\leqslant \{1+o(1)\}\widehat{g}(0)X\int_{\frac{1}{2}}^{\vartheta}\frac{2}{\gamma(\theta)}\ud\theta,
\]
where $\gamma(\theta)$ is given by \eqref{eq:gamma(theta)}.
A numerical calculation shows that
\begin{align*}\int_{\frac{1}{2}}^{\frac{64}{97}}\frac{2\cdot62}{91-89\theta}\ud\theta+\int_\frac{64}{97}^\frac{32}{41}\frac{2\cdot60}{86-83\theta}\ud\theta+\int_\frac{32}{41}^\vartheta\frac{2\cdot14}{19-18\theta}\ud\theta<\frac{3}{2}\end{align*}
with $\vartheta=0.847.$ This, together with (\ref{eq:H(X)}), (\ref{eq:H(X)2}) and (\ref{eq:H1(X)bound}), implies 
\begin{align}\label{eq:H3(X)}
H_3(X)\gg X\end{align}
for such $\vartheta,$
provided that
\begin{align}\label{eq:H4(X)}
H_4(X)=o(X).\end{align}
We then conclude from (\ref{eq:H3(X)}) that $P^+(p^2+1)>p^{0.847}$ for infinitely many primes $p$, proving Theorem \ref{thm:P^+}.

It remains to prove (\ref{eq:H4(X)}) and this can be concluded from the square sieve of Heath-Brown as follows. Note, for any fixed $\ell\geqslant1,$ that
\begin{align*}Q_\ell(X)\ll\Big(\frac{X}{\ell}+1\Big)\rho(\ell),\end{align*}
from which we may derive trivially that
\begin{align*}\sum_{k\geqslant3}\sum_{X^\flat<p^k\ll X^2}Q_{p^k}(X)\log p
&\ll X^\varepsilon\sum_{3\leqslant k\leqslant3\log X}\sum_{X^\flat<p^k\ll X^2}\Big(\frac{X}{p^k}+1\Big)
\\\noalign{\vskip -0,5mm}
&\ll X^{1/2+\varepsilon}+ X^\varepsilon\sum_{3\leqslant k\leqslant3\log X}\sum_{p\ll X^{2/3}}1
\\\noalign{\vskip -1,5mm}
&\ll X^{2/3+\varepsilon}.
\end{align*}
It remains to show that
\begin{align*}
\sum_{\sqrt{X^\flat}<p\ll X}Q_{p^2}(X)\log p=o(X).\end{align*}
In fact, we shall prove the following slightly stronger estimate
\begin{align*}
\cN(X)
:= \sum_{\ell\sim L} \sum_{\substack{n\sim X\\n^2+1\equiv0\!\bmod{\ell^2}}}1\ll X^{1-\varepsilon}\end{align*}
for all $\sqrt{X^\flat}\ll L\ll X$. For $\sqrt{X^\flat}\ll L\ll X^{1-2\varepsilon}$, the above argument also applies. We only consider the remaining case $X^{1-2\varepsilon}\ll L\ll X$. 
In fact, the bound for $\cN(X)$ was already obtained in \cite{Da} appealing to the following
square sieve of Heath-Brown \cite{HB}. For the completeness of arguments, we present the proof as quickly as possible.

\begin{lemma}[Square sieve]
Let $\xi:\mathbf{N}\rightarrow\bR_{\geqslant0}$ be an arbitrary function with 
$\sum_{n\in\mathbf{N}}\xi(n)<\infty.$ 
Suppose $\mathfrak{P}$ is a set of $P$ prime numbers and $\xi$ vanishes if $n=0$ 
or $n\geqslant \mathrm{e}^P,$ then we have
\begin{align*}
\sum_{n\in\mathbf{N}} \xi(n^2)
\ll \frac{1}{P}\sum_{n\in\mathbf{N}}\xi(n)
+ \frac{1}{P^2}\mathop{\sum\sum}_{p\neq q\in\mathfrak{P}}
\bigg|\sum_{n\in\mathbf{N}}\xi(n)\bigg(\frac{n}{pq}\bigg)\bigg|,
\end{align*}
where $(\frac{\cdot}{pq})$ denotes the Jacobi symbol $\bmod{pq}.$
\end{lemma}

We now introduce a set $\mathfrak{P}$ consisting of $P$ prime numbers, where $P$ is a large number to be specialized later.
For $n^2+1\equiv0\bmod{\ell^2},$ we can write $n^2+1=m\ell^2$ for some $m\in (M,2M]$ 
with $M\asymp X^2L^{-2}\ll X^{4\varepsilon}.$ Thus, we may apply the square sieve to the sequence 
$\{m\ell^2-1\}_{\ell\sim L, \, m\sim M}$. More precisely, we have
\begin{align*}
\cN(X)
& \ll \frac{1}{P}\sum_{\ell\sim L}\sum_{m\sim M}1
+ \frac{1}{P^2} \mathop{\sum_{p\in\mathfrak{P}}\sum_{q\in\mathfrak{P}}}_{p\neq q}
\bigg|\sum_{\ell\sim L}\sum_{m\sim M} \bigg(\frac{m\ell^2-1}{pq}\bigg)\bigg|
\\
& \ll \frac{LM}{P} + \frac{1}{P^2} \mathop{\sum_{p\in\mathfrak{P}}\sum_{q\in\mathfrak{P}}}_{p\neq q}
\sum_{m\sim M} \bigg|\sum_{\ell\sim L} \bigg(\frac{m\ell^2-1}{pq}\bigg)\bigg|.
\end{align*}
By completing method and Weil's bound for complete character sums, the innermost sum over $\ell$ is bounded by $(pq)^{1/2+\varepsilon}(m,pq)$. Therefore, for $X^{1-2\varepsilon}\ll L\ll X$,
\begin{align*}
\cN(X)&\ll LMP^{-1}+PX^{5\varepsilon}\ll X^{1-\varepsilon}\end{align*}
on taking $P=X^{5\varepsilon}.$ This completes the proof of (\ref{eq:H4(X)}), thus that of Theorem \ref{thm:P^+}.

\bigskip

\section{Improving a result of Dartyge} \label{sec:Dartyge}

The proof of Theorem \ref{thm:Dartyge+} also follows from the Chebyshev-Hooley method.
Before starting the proof, we would like to recall the counting function of sifted numbers.
Write $\Phi(X,z)=\Phi(X,z;1),$ so that
\begin{align*}
\Phi(X,z) := \sum_{ p\mid n\Rightarrow p>z}g\Big(\frac{n}{X}\Big).
\end{align*}
A classical result, $u:=\log X/\log z,$
\begin{align*}
\Phi(X,z)
= \widehat{g}(0)\frac{Xw(u)-z}{\log z}+O\bigg(\frac{X}{(\log z)^2}\bigg),
\end{align*}
where $w(u)$ is the Buchstab function defined recursively by
\begin{align*}
\begin{cases}
uw(u)=1              & (1\leqslant u\leqslant2),
\\
(uw(u))'=w(u-1)   & (u>2).
\end{cases}
\end{align*}

Suppose now  $1<u\leqslant13$, and we would like to examine the sum
\begin{align*}
H(u,X) := \sum_{p\mid n\Rightarrow p>X^{1/u}}g\Big(\frac{n}{X}\Big)\log(n^2+1).
\end{align*}
On one hand, we have
\begin{align*}
H(u,X)
& = \{2\log X+O(1)\} \Phi(X,X^{1/u})
= 2uw(u)\widehat{g}(0)X \{1+o(1)\}.
\end{align*}
On the other hand, 
the relation \eqref{eq:log=Lambda} allows us to write
\begin{align*}
H(u,X) 
= \sum_{\ell\ll X^2} \Lambda(\ell)
\sum_{\substack{p\mid n\Rightarrow p>X^{1/u}\\ n^2+1\equiv0\!\bmod\ell}} g\Big(\frac{n}{X}\Big)
= \sum_{\ell\ll X^2} \Lambda(\ell) Q_\ell(X; u)
\end{align*}
and split the sum over $\ell$ following the manner in \eqref{eq:H(X)2}.
By virtue of Lemmas \ref{lm:BV-Wolke} and \ref{lm:Brun-Titchmarsh-generalized},
it suffices to find the smallest $u>1$ such that
\begin{align*}
& \int_{\frac{1}{2}}^{\frac{16}{17}}F(u\gamma(\theta))\ud\theta
+ \int_{\frac{16}{17}}^{\theta_0}F(u(1-\theta))\ud\theta
+ \frac{u}{\mathrm{e}^{\gamma}}
\int_{\theta_0}^1\frac{\theta \ud\theta}{\sigma_2((\frac{2}{3}-\frac{\theta}{2})u)}
<\frac{3}{2}\mathrm{e}^{\gamma}w(u),
\end{align*}
where $\theta_0=0.9926$ as chosen in \cite[Section 9]{Da}, $\gamma(\theta)$ is given by \eqref{eq:gamma(theta)}, the second integral comes from the classical Brun-Titchmarsh theorem of van Lint-Richert and $1/\sigma_2(s)$ appears in the Selberg sieve of dimension 2, 
which is equal to $8\mathrm{e}^{2\gamma}s^{-2}$ if $0<s\leqslant2$.
One may check with the help of Mathematica 9 that $u=11.2$ works.

\bigskip

\section{Proof of Theorem \ref{thm:P_4}} \label{sec:proof-P_4}

\subsection{Preparation for sifting}
We first state some convention to sift the specialized sequence
\[
\sA := \big\{\tfrac{1}{2}(p^2+1):X<p\leqslant2X\big\},
\]
although most arguments are suited for all general non-negative sequences.

Define the {\it smoothed} sifting function
\begin{align*}
S(\sA,z) := \sum_{(\frac{1}{2}(p^2+1),P(z))=1}g\Big(\frac{p}{X}\Big),
\end{align*}
where, for $z>3$,
\[P(z) := \prod_{2<p<z}p.\]
For squarefree $d$, we consider subsequence 
\[
\sA_d := \big\{\tfrac{1}{2}(p^2+1):X<p\leqslant2X \;\text{and}\;\, p^2+1\equiv0\bmod d\big\}.
\]
Its sifting function is defined by
\begin{align*}
S(\sA_d,z)
= \sum_{\substack{(\frac{1}{2}(p^2+1),P(z))=1\\p^2+1\equiv 0\!\bmod d}}g\Big(\frac{p}{X}\Big).
\end{align*}

Recall congruence sum $Q_d(X)$, defined by \eqref{eq:Q_l(X)}. 
For $d\leqslant X^\flat$, the Bombieri-Vinogradov theroem (see Lemma \ref{lm:BV}) yields that $Q_d(X)$ can be approximated on average by $Q_1(X)\rho(d)/\varphi(d)$.

The following lemma then characterizes the dimension of sieves, see \cite[Proposition 10.1]{DH} for instance.

\begin{lemma}\label{lm:sieve-dim} 
For $z>3,$ we have  
\begin{equation}\label{cond:dimension}
\sum_{p\leqslant z}\frac{\rho(p)}{\varphi(p)}\log p=\log z+O(1)
\end{equation}
and
\begin{equation}\label{def:Vz}
V(z) := \prod_{2<p\leqslant z} \bigg(1-\frac{\rho(p)}{\varphi(p)}\bigg)
= \{1+o(1)\} \frac{\mathrm{e}^{-\gamma} \fc}{\log z}
\end{equation}
with
\begin{align*}
\fc := 2 \prod_{p>2} \bigg(1-\frac{\rho(p)}{\varphi(p)}\bigg)\bigg(1-\frac{1}{p}\bigg)^{-1}.
\end{align*}
\end{lemma}

\subsection{A weighted sieve}
We now introduce the weighted sieve of Richert following the manner in \cite{Ri}.
Let $0<\alpha<\beta$ be some constants to be specialized later. Put
\begin{align*}
z := X^\alpha,
\qquad
y := X^\beta.
\qquad
\eta := r+1-2/\beta.
\end{align*}

Consider the weighted sum
\begin{align*}
\Psi(\alpha,\beta;\eta)
:= \sum_{(\frac{1}{2}(p^2+1),P(z))=1} g\Big(\frac{p}{X}\Big)
\Big(1-\frac{1}{\eta}\sum_{\substack{q\mid (p^2+1)\\ z\leqslant q<y}}w_q\Big)
\end{align*}
with $w_q := 1-(\log q)/\log y$.
Note that $q$ is a prime variable. 

The following lemma is taken from \cite{Ri} with slight modification in notation.
\begin{lemma}\label{lm:weightedsieve}
Suppose that for a given $r\geqslant1$, there exist constants $\alpha,\beta$ with
\[
0<\alpha<\beta<1,\qquad
\beta>2/(r+1)
\]
such that
\begin{align}\label{eq:Phi-lowerbound}
\Psi(\alpha,\beta;\eta)\gg X(\log X)^{-2}.
\end{align}
Then we have
\begin{align*}
\big|\big\{X<p\leqslant2X \,:\, \tfrac{1}{2}(p^2+1)=P_r\big\}\big|\gg X(\log X)^{-2}.
\end{align*}
\end{lemma}

Theorem \ref{thm:P_4}, in the case of such special polynomial, follows from suitable choices for $\alpha,\beta$ such that 
\eqref{eq:Phi-lowerbound} holds with $r=4$. To do so, we will seek the lower bound for $\Phi(\alpha,\beta;\eta)$ starting from the following 
expression
\begin{align}\label{eq:Phi(alpha,beta)}
\Psi(\alpha,\beta;\eta)= S(\sA,z) - \eta^{-1} \sum_{z\leqslant q<y}w_qS(\sA_q,z).
\end{align}

\subsection{Sieve estimates}

The lower bound for $S(\sA,z)$ follows from a routine application of lower-bound sieve 
and Lemma \ref{lm:BV}.

\begin{lemma}\label{lm:lower}
For $X\to+\infty$, we have
\begin{align*}
S(\sA,z)\geqslant f\Big(\frac{1}{2\alpha}\Big)\widehat{g}(0)V(z)\frac{X}{\log X} \{1+o(1)\}
,\end{align*}
where $f$ is defined by the system $\eqref{eq:sieve-F}$
and $V(z)$ is given by \eqref{def:Vz}.
\end{lemma}

An upper bound for $S(\sA_q,z)$ with small primes $q$ follows from the upper-bound sieve and 
Lemma \ref{lm:BV}. 
For larger primes $q$, especially while going beyond the Bombieri-Vinogradov theorem, 
we will appeal to Lemma \ref{lm:error-average} combining with a composition of two linear sieves. 
To do so, we introduce a parameter $\delta\in (\alpha, \beta)$ to be optimized.

\begin{lemma}\label{lm:upper1}
For $X\to+\infty$, we have
\begin{align*}
\sum_{z\leqslant q<X^\delta}w_qS(\sA_q,z)
\leqslant c_1 \widehat{g}(0)V(z)\frac{X}{\log X} \{1+o(1)\},
\end{align*}
where $V(z)$ is given by \eqref{def:Vz} and 
\begin{align}\label{eq:c1}
c_1 := \int_\alpha^\delta\Big(\frac{1}{\theta}-\frac{1}{\beta}\Big)F\Big(\frac{1-2\theta}{2\alpha}\Big)\ud\theta.\end{align}
\end{lemma}

\proof
In fact, an upper-bound sieve of Rosser-Iwaniec yields
\begin{align*}
S(\sA_q,z)
& \leqslant \frac{X}{\log X}\widehat{g}(0)V(z)\frac{\rho(q)}{\varphi(q)}
\bigg\{F\bigg(\frac{\log(X^\flat/q)}{\log z}\bigg)+O\bigg(\frac{1}{(\log X)^{1/6}}\bigg)\bigg\}
\\
& \qquad 
+ O\bigg(\sum_{d\leqslant X^\flat/q}\bigg|Q_{qd}(X)-\frac{\rho(qd)}{\varphi(qd)}Q_1(X)\bigg|\bigg),\end{align*}
where $F$ is defined recursively by \eqref{eq:sieve-F}.
Summing over $q$ against the weight $w_q$ with integration by parts thanks to \eqref{cond:dimension}
and controlling the error term by virtue of Lemma \ref{lm:BV},
we get the required result.
\endproof

We now turn to consider $S(\sA_q,z)$ for $q\geqslant X^\delta.$ We will appeal to
Lemma \ref{lm:error-average} instead of Lemma \ref{lm:BV}.
\begin{lemma}\label{lm:upper2}
Let $\beta<0.68$. For $X\to+\infty$, we have
\begin{align*}
\sum_{X^\delta\leqslant q<X^\beta}w_qS(\sA_q,z)
\leqslant c_2\widehat{g}(0)V(z)\frac{X}{\log X} \{1+o(1)\},
\end{align*}
where
\begin{align}\label{eq:c2}
c_2
& :=\mathrm{e}^{\gamma} \alpha \int_\delta^\beta\Big(\frac{1}{\theta}-\frac{1}{\beta}\Big)
\frac{88288\, \ud\theta}{8281 - 16198 \theta + 7921 \theta^2}\cdot
\end{align}
\end{lemma}

\proof
We would like to make initial estimates for each $S(\sA_q,z)$ from above 
by the composition of two upper-bound linear sieves (see Theorem \ref{thm:sieve-composition} below), 
and then all parameters will be optimized after summing over $q$.

For $Q=X^\theta$ with $\delta\le \theta\le \beta$,
define
\begin{align*}
\sC(Q):=\sum_{Q<q\le 2Q} w_qS(\sA_q,z).
\end{align*}
Let $\lambda_1,\lambda_2$ be two upper-bound linear sieves, of level $D_1,D_2$, 
so that $0\le 1*\mu\le 1*\lambda_i$ for $i=1,2$. 
We thus have, with $P_*=\prod_{p\le \sqrt{X}} p$
and the notation \eqref{def:AdXell},
\begin{align*}
S(\sA_q,z)
& = \sum_{d_1\mid P(z)} \mu(d_1) \sum_{\substack{d_2\mid P_*\\(d_2,qd_1)=1}} \mu(d_2) A_{d_2}(X; 2d_1q)
\\
& \leqslant \sum_{\substack{d_1\leqslant D_1\\d_1|P(z)}}\sum_{\substack{d_2\leqslant D_2\\(d_2,qd_1)=1}}\lambda_1(d_1)\lambda_2(d_2)A_{d_2}(X; 2d_1q).
\end{align*}
The conditions that $d_1|P(z)$ and $q\geqslant z$ implies that $(d_1,q)=1.$
Replacing $A_{d_2}(X; 2d_1q)$ by $\widehat{g}(0)\rho(d_1)\rho(q)(d_1d_2q)^{-1}X + r_{d_2}(X;2qd_1)$
(cf. \eqref{eq:error-r_d}) and inserting the obtained inequality into the definition of $\mathscr{C}(Q)$,
we find that
\begin{align}\label{eq:sC(Q)}
\sC(Q)&\leqslant\sM(Q)+\sE(Q)\end{align}
with
\begin{align*}
\sM(Q)
& := \widehat{g}(0)X\sum_{Q<q\le 2Q} w_q \frac{\rho(q)}{q}
\sum_{\substack{d_1\leqslant D_1\\d_1|P(z)}}\frac{\lambda_1(d_1)\rho(d_1)}{d_1}
\sum_{\substack{d_2\leqslant D_2\\(d_2,qd_1)=1}} \frac{\lambda_2(d_2)}{d_2},
\\
\sE(Q)
& := \sum_{Q<q\le 2Q} w_q\sum_{\substack{d_1\leqslant D_1\\d_1|P(z)}}\sum_{\substack{d_2\leqslant D_2\\(d_2,qd_1)=1}}\lambda_1(d_1)\lambda_2(d_2)r_{d_2}(X;2qd_1).\end{align*}

Thanks to Iwaniec \cite{Iw2}, the sieve weights $\lambda_1,\lambda_2$ can be chosen to be finite linear combinations of some functions that are well-factorable of degree $J$ for any fixed large $J\geqslant2$, so that we are in a position to apply Lemma
\ref{lm:error-average}, getting
\begin{align}\label{eq:sE(Q)}
\sE(Q)&\ll X^{1-\delta}\end{align}
for some $\delta>0$, provided that
\begin{align}\label{eq:constraints}
\gamma_2<\tfrac{91}{62}-\tfrac{89}{62}(\gamma_1+\theta),
\qquad 
\tfrac{1}{2}\leqslant\gamma_1+\theta<\tfrac{112}{131}\end{align}
with
\[D_1=X^{\gamma_1},
\qquad 
D_2=X^{\gamma_2},
\qquad 
Q=X^\theta.\]

The upper bound for $\sC(Q)$ will be established by evaluating the main term $\sM(Q)$. This is in fact a composition of two linear upper-bound sieves, and we appeal to a reduction of Friedlander and Iwaniec, see Theorem \ref{thm:sieve-composition} in the appendix, in which we should take
\begin{align*}
g_1(d) = \begin{cases} 
\rho(d)/d & (2\nmid d)
\\
0             & (2\mid d)
\end{cases}
\qquad\text{and}\qquad
g_2(d) = \begin{cases} 
1/d & (q\nmid d)
\\
0    & (q\mid d)
\end{cases}.
\end{align*}
It is easy to check that both of $g_1,g_2$ satisfy the restriction \eqref{eq:g1g2}.
Hence it follows 
\begin{align}\label{eq:sM(Q)}
\sM(Q)
\leqslant \{1+o(1)\} \widehat{g}(0)\frac{4\mathrm{e}^{\gamma}\alpha}{\fc \gamma_1\gamma_2}
\frac{XV(z)}{\log X}
\sum_{Q<q\le 2Q} w_q \frac{\rho(q)}{q} H_q,
\end{align}
where $\fc$ is defined in Lemma \ref{lm:sieve-dim} and 
$H_q$ is corresponds to $H$ in Theorem \ref{thm:sieve-composition} upon the above choices for $g_1,g_2$, i.e.,
\begin{align*}
H_q
& = 2\bigg(1-\frac{\rho(q)}{q}\bigg) \bigg(1-\frac{1}{q}\bigg)^{-2}
\prod_{p\nmid2q}\bigg(1-\frac{1+\rho(p)}{p}\bigg)\bigg(1-\frac{1}{p}\bigg)^{-2}
\\\noalign{\vskip -0,5mm}
& =2\bigg(1-\frac{\rho(q)}{q}\bigg)\bigg(1-\frac{1+\rho(q)}{q}\bigg)^{-1}\prod_{p>2}\bigg(1-\frac{1+\rho(p)}{p}\bigg)\bigg(1-\frac{1}{p}\bigg)^{-2},\end{align*}
which can be reduced to 
\begin{align*}
H_q=\fc\cdot\bigg(1-\frac{\rho(q)}{q}\bigg)\bigg(1-\frac{1+\rho(q)}{q}\bigg)^{-1}
\end{align*}
since $\rho(p)=2$ for $p\equiv 1\bmod 4$ and $\rho(p)=0$ for $p\equiv 3\bmod 4$.

It suffices to maximize $\gamma_1\gamma_2$, in terms of $\theta$,
subject to the constraints in \eqref{eq:constraints}.
With the help of Mathematica 9, we have
\begin{align*}
\gamma_1\gamma_2<\tfrac{1}{22072}(8281 - 16198 \theta + 7921 \theta^2)
\quad\mathrm{for}\quad
\theta\in(0, \tfrac{8015}{11659}).
\end{align*}
One may see why we require $\beta<0.68<\tfrac{8015}{11659}$ in Lemma \ref{lm:upper2}.

Collecting all $Q$, Lemma \ref{lm:upper2} then follows from \eqref{eq:sC(Q)}, \eqref{eq:sE(Q)}, \eqref{eq:sM(Q)} and partial summation.
\endproof

\subsection{Conclusion of Theorem \ref{thm:P_4}}
We conclude from Lemmas \ref{lm:lower}-\ref{lm:upper2}
that
\begin{align}\label{eq:Phi2}
\Psi(\alpha,\beta;\eta)\geqslant C\widehat{g}(0)
V(z)\frac{X}{\log X} \{1+o(1)\}
\end{align}
with
\begin{align*}
C
& := f\Big(\frac{1}{2\alpha}\Big)-\frac{1}{\eta}
\int_\alpha^\delta\Big(\frac{1}{\theta}-\frac{1}{\beta}\Big)F\Big(\frac{1-2\theta}{2\alpha}\Big)\ud\theta
\\
& \quad
- \frac{\mathrm{e}^{\gamma} \alpha}{\eta}
\int_\delta^\beta\Big(\frac{1}{\theta}-\frac{1}{\beta}\Big) 
\frac{88288\, \ud\theta}{8281 - 16198 \theta + 7921 \theta^2}\cdot
\end{align*}

For $\theta\geqslant\frac{1}{2}-3\alpha,$ we find
\[
F\Big(\frac{1-2\theta}{2\alpha}\Big)
= \frac{4\mathrm{e}^{\gamma}\alpha}{1-2\theta}\cdot
\]
Thus, $\delta\approx0.44$ is chosen to be the root of the equation
\begin{align*}
\frac{1}{1-2\delta}
= \frac{22072}{8281 - 16198\delta + 7921 \delta^2}\cdot
\end{align*}

We choose $\alpha=\frac{1}{12}$, so that
\begin{align*}
\frac{1-2\theta}{2\alpha}
= 6(1-2\theta)
\in \begin{cases}
[1,3] & \text{if $\,\theta\in[\frac{1}{4},\delta]$},
\\\noalign{\vskip 1mm}
[3,5] & \text{if $\,\theta\in[\frac{1}{12},\frac{1}{4}]$}.
\end{cases}
\end{align*}
Note that
\begin{align*}
F(s) = \begin{cases}
\displaystyle\frac{2\mathrm{e}^{\gamma}}{s} & (1\le s\le 3),
\\\noalign{\vskip 1mm}
\displaystyle\frac{2\mathrm{e}^{\gamma}}{s}\bigg(1+\int_2^{s-1}\frac{\log(t-1)}{t}\ud t\bigg) & (3\le s\le 5).
\end{cases}
\end{align*}
This will lead to the following more explicit expression for $C$
\begin{align*}
C
& = f(6)
- \frac{\mathrm{e}^{\gamma}}{3\eta}
\bigg\{\int_{\frac{1}{12}}^\frac{1}{4}\Big(\frac{1}{\theta}-\frac{1}{\beta}\Big)\Big(1+\int_2^{6(1-2\theta)-1}\frac{\log(t-1)}{t}\ud t\Big)\frac{\ud\theta}{1-2\theta}
\\
& \quad
+ \int_{\frac{1}{4}}^{0.44}\Big(\frac{1}{\theta}-\frac{1}{\beta}\Big)\frac{\ud\theta}{1-2\theta}
+ \int_{0.44}^{\beta} \Big(\frac{1}{\theta}-\frac{1}{\beta}\Big)
\frac{22072\, \ud\theta}{8281 - 16198 \theta + 7921 \theta^2}
\bigg\}.
\end{align*}

Taking $r=4$ and $\beta=0.622$, we find $C\approx0.0568.$
This establishes Theorem \ref{thm:P_4} in view of Lemma \ref{lm:weightedsieve}.

\begin{figure}
\includegraphics[width=10cm]{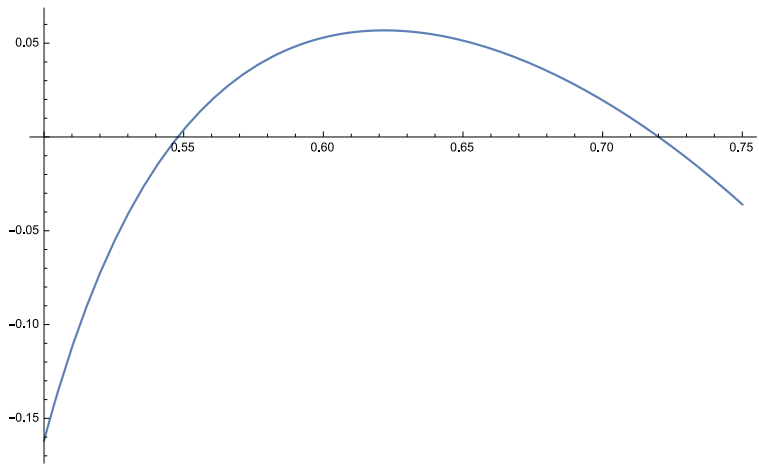}
\caption{Graph for $C=C(\beta)$ with $r=4$}
\end{figure}

\vskip 8mm

\appendix

\section{Composition of two linear sieves}
In this appendix, we formulate a reduction of Friedlander and Iwaniec on the composition of two sieves. In particular, we focus on linear sieves. The original statement with proof can be found in \cite[Appendix A]{FI1} or \cite[Section 5.10]{FI2}.

\begin{theorem}\label{thm:sieve-composition}
Let $(\lambda_1),$ $(\lambda_2)$ be two upper-bound linear sieves of levels $D_1,D_2,$ respectively, 
and let $g_1,g_2$ be density functions satisfying linear sieve conditions and 
\begin{align}\label{eq:g1g2}
0\leqslant g_1(p),g_2(p)\leqslant\tfrac{1}{2}
\end{align}
for each prime $p$. Then
$$
\bigg|\mathop{\sum\sum}_{(d_1,d_2)=1}\lambda_1(d_1)\lambda_2(d_2)g_1(d_1)g_2(d_2)\bigg|
\leqslant \frac{4H\{1+o(1)\}}{\log D_1\log D_2}
$$
with
$$
H := \prod_{p}(1-g_1(p)-g_2(p))(1-1/p)^{-2}.
$$
\end{theorem}

One can see that we have an extra condition \eqref{eq:g1g2} compared to the original version of Friedlander and Iwaniec, for whom the constant $H$ should be replaced by the following larger one
$$
\prod_{p}(1-g_1(p)-g_2(p)+2g_1(p)g_2(p))(1-1/p)^{-2}.
$$
The modification here is due to avoiding the use of the trivial inequality
\[\left|1-\frac{g_1g_2}{(1-g_1)(1-g_2)}\right|\leqslant1+\frac{g_1g_2}{(1-g_1)(1-g_2)}\]
at primes.

\end{document}